\documentclass{amsart}

\newtheorem{teo}{Theorem}[section]
\newtheorem{prop}[teo]{Proposition}
\newtheorem{cor}[teo]{Corollary}

\newtheorem{claim}[teo]{Claim}

\theoremstyle{definition}
\newtheorem{definicao}[teo]{Definition}

\theoremstyle{remark}
\newtheorem{rem}[teo]{Remark}

\numberwithin{equation}{section}

\newcommand{\dist}{\ensuremath{\mathrm{dist} }}
\newcommand{\grad}{\ensuremath{\mathrm{grad}\ }}

\newcommand{\tub}{\ensuremath{\mathrm{Tub} }}

\newcommand{\F}{\ensuremath{\mathcal{F}}}
\newcommand{\singularF}{\ensuremath{\mathcal{X}_{F}}}

\newcommand{\eps}{\varepsilon}

\begin{document}

\title{Generalizations of Isoparametric Foliations}

\author{Marcos M. Alexandrino} 
 

\address{Marcos M. Alexandrino\\Departamento de Matem\'{a}tica, PUC-Rio, 
 Rua Marqu\^{e}s de S\~{a}o Vicente, 225, 
  22453-900, Rio de Janeiro, Brazil}
\email{:malex@mat.puc-rio.br}
\email{marcosmalex@yahoo.de}

\subjclass{Primary 53C12, Secondary 57R30}

\date{October, 2004.}

\keywords{Singular Riemannian foliations, isoparametric foliations,  orbifolds, equifocal submanifolds,  polar actions.}


\maketitle

\section{Introduction}

Isoparametric submanifolds and hypersurfaces in space forms are geometric objects that have been studied since \'{E}. Cartan.  

Another important class of geometric objects  is the orbits of polar actions on a Riemannian manifold, e.g., the orbits of the adjoint action of a compact Lie group on itself.

These two classes of submanifolds share  some common properties. For example,  they are leaves of  \emph{singular Riemannian foliations with section} (s.r.f.s. for short). 

A singular foliation on a complete Riemannian manifold is said to be Riemannian if every geodesic that is perpendicular at one point to a given leaf remains perpendicular to every leaf it meets. Moreover, the singular foliation admits sections if, for each regular point $p,$ there is a totally geodesic complete immersed submanifold through $p$ that meets each leaf orthogonally and whose dimension is the codimension of the regular leaves.

The purpose of this paper is to review some results of the theory of s.r.f.s., introduced in \cite{Alex3} and  developed in \cite{Alex1},\cite{Alex2} and \cite{AlexToeben}.
This paper is organized as follows.
In Section 2  we recall the definitions of isoparametric submanifolds, polar actions and equifocal submanifolds (the last one  introduced by Terng and Thorbergsson \cite{TTh1} in order to generalize isoparametric submanifolds in symmetric spaces). Then  s.r.f.s. is presented as the natural candidate to generalize all these concepts.
In Section 3 another possible way to generalize isoparametric foliations is presented. We consider the set of  preimages of isoparametric and transnormal maps as well  some results about these maps. Then  we recall Theorem \ref{TransormalMapTeo}
, which states that the preimages of regular values are leaves of a s.r.f.s.  
In Section 4  some results  about s.r.f.s. on Riemannian manifolds are listed. The main results of this section  are Theorem \ref{frss-eh-equifocal}
 and Theorem \ref{sliceteorema}.
In Section 5  a joint work with Dirk T\"{o}ben is briefly discussed. Inspired by simple examples and some results of Lie Group theory, we ask if there exists an exceptional leaf of a s.r.f.s. on a simply connected space and if there exist fundamental domains in the sections. These questions are answered in Theorem \ref{teo-excepcional-leaf} and Theorem \ref{teo-fundamental-domain}.
Finally, in Section 6, I give  an alternative proof of T\"{o}ben's Theorem (see Theorem \ref{teo-equical-srfs}), who gave in  
 \cite{Toeben2} a necessary and sufficient condition for an equifocal submanifold  to induce a s.r.f.s.  

I would like to take the opportunity to thank  Paolo Piccione, Claudio Gorodski and the other organizers of XIII Escola de Geometria Diferencial for their invitation and hospitality in USP.
 I also thank Julio Rebelo for useful suggestions.

\section{The natural Candidate to Generalize Isoparametric Foliations}

For the sake of motivation, let us start by recalling a classical result from the Lie Group theory (see e.g. the book of Duistermaat and Kolk \cite{Duistermaat}). 

\begin{teo}[Maximal Torus Theorem]
Let $G$ be a compact Lie group with a bi-invariant metric and consider the action of $G$ on  itself by conjugation (the so called \textbf{adjoint action}). Then
\begin{enumerate}
\item[1)] There exist maximal tori.
\item[2)] The orbits of the adjoint action meet each maximal torus orthogonally.
\item[3)] The intersection of the singular orbits with  a fixed torus $T$ is a finite union of totally geodesic hypersurfaces (the walls).
\item[4)] Let $x$ be a point of  $T$ and $G(x)$ its orbit. Then $G(x)\cap T$ is invariant by reflections in the walls.
\end{enumerate}
\end{teo}

The adjoint action is a typical example of a polar action, whose definition is recalled below.

\begin{definicao}
An isometric action of a compact Lie group $G$ on a complete Riemannian manifold $M$ is called \emph{polar} if there exists a complete immersed submanifold $\Sigma$ of $M$ that meets all $G$-orbits orthogonally. Such $\Sigma$ is called \emph{section}. 
If the sections are flat the action is called \emph{hyperpolar}. 
\end{definicao}

Apart from the classical results about the adjoint action, there are nowadays several works about polar actions, in particular about polar actions on symmetric spaces. Among these works, we would like to cite Kollross \cite{Kollross}, where a classification of hyperpolar actions on irreducible symmetric spaces of compact type is given. Examples of polar actions that are not hyperpolar and a classification of polar action on compact rank one symmetric spaces can be found in the paper of  F. Podest\`{a} and G. Thorbergsson \cite{PodestaThorbergsson}.


The orbits of a polar action on an Euclidean space are examples of isoparametric submanifolds, whose concept we recall below.

Let $L$ be an immersed submanifold of a Riemannian manifold $M.$  A section $\xi$ of  the normal bundle $\nu (L)$ is said  to be a \emph{parallel normal field} along $L$ if  $\nabla^{\nu}\xi\equiv 0,$ where $\nabla^{\nu}$ is the normal connection.  $L$  is said to have globally flat normal bundle, if the holonomy  of the normal bundle  $\nu (L)$ is trivial, i.e., if any normal vector can be extended to a globally defined parallel normal field. 

\begin{definicao}
A submanifold $F$ of a space form $M(k)$ is called \emph{isoparametric} if its normal bundle is flat and if the principal curvatures along any parallel normal vector field are constant.
\end{definicao} 

Isoparametric hypersurfaces in space forms have been studied since \'{E}. Cartan \cite{Cartan1,Cartan2,Cartan3,Cartan4}. In the eighties C. E. Harle\cite{Harle}, S. Carter and A. West \cite{CarterWest1, CarterWest2} and C-L.Terng \cite{Terng} independently introduced the concept of isoparametric submanifold in a space form. Terng observed in her work the similarity between  isoparametric submanifolds and the  polar orbits in Euclidean spaces (e.g. adjoint orbits in the Lie algebra). In particular, she associated Coxeter groups  to isoparametric submanifolds.  

An interesting result that relate orbits of polar action to isoparametric submanifold is due to G. Thorbergsson.  In \cite{Th2}  Thorbergsson proved that a  compact, irreducible isoparametric submanifold $L^{n}$ in $\textbf{R}^{n+k}$ is homogeneous if $k\geq 3$ and $L^{n}$ does not lie in any affine hyperplane of $\textbf{R}^{n+k}.$ Nevertheless, there is a lot of non homogeneous examples of codimension 2 in Euclidean spaces (isoparametric hypersurfaces in spheres).

A comprehensive description of Terng's work about isoparametric submanifolds on space forms and Hilbert spaces can be found in the book of Palais and Terng \cite{PTlivro}.
The history of isoparametric hypersurfaces and submanifolds and their generalizations until  2000 can be found in the survey  of  Thorbergsson \cite{Th}.

%

In \cite{TTh1} C.L. Terng and G. Thorbergsson introduced the concept of equifocal submanifold with  flat sections in symmetric space in order to generalize the definition of isoparametric submanifold in Euclidean space. Now we review the definition of equifocal submanifolds in a Riemannian manifold. 

\begin{definicao} 
A connected immersed submanifold  $L$  of a complete Riemannian manifold $M$ is called \textbf{equifocal} if
\begin{enumerate}
\item[0)] the normal bundle $\nu(L)$ is globally flat,
\item[1)] $L$ has sections, i.e., for all $p\in L$ there exists a complete, immersed, totally geodesic submanifold $\Sigma$  such that $\nu_{p}(L)=T_{p}\Sigma.$
\item[2)] For each  parallel normal field $\xi$  along $L,$  the derivative of  the map    $\eta_{\xi}:L\rightarrow M,$ defined as $\eta_{\xi}(x):=\exp_{x}(\xi),$ has constant rank.
\end{enumerate}
\end{definicao}
A connected immersed submanifold $L$ is called \emph{local equifocal} if, for each $p\in L,$ there exists a neighborhood $U\subset L$ of $p$ in $L$ such that $U$ is an equifocal submanifold.

According to Terng, Palais and Thorbergsson (see \cite{PTlivro} and \cite{TTh1}), the equifocal submanifolds with flat section on a compact symmetric space, isoparametric submanifolds on an Euclidean space and orbits of polar actions on a Riemannian manifold share  a common property. They are leaves of  singular Riemannian foliations with sections, whose definition is given  below. 

At first  we must recall the concept of singular Riemannian foliation (see the book of Molino \cite{Molino}).

\begin{definicao}
 A partition $\F$ of a complete Riemannian manifold $M$ by connected immersed submanifolds (the \emph{leaves}) is called a \emph{singular foliation} of $M$ if it verifies condition (1) and {singular Riemannian foliation} if it verifies condition (1) and (2):

\begin{enumerate}
\item $\F$ is \emph{singular}
i.e., the module $\singularF$ of smooth vector fields on $M$ that are tangent at each point to the corresponding leaf acts transitively on each leaf. In other words, for each leaf $L$ and each $v\in TL$ with foot point $p$ there is a $X\in \singularF$ with $X(p)=v$.
\item  The partition is \emph{transnormal}, i.e., every geodesic that is perpendicular at one point to a leaf remains perpendicular to every leaf it meets.
\end{enumerate}
\end{definicao}

Let $\F$ be a singular Riemannian foliation on an complete Riemannian manifold $M.$  A point $p\in M$ is called \emph{regular} if the dimension of the leaf $ L_{p}$ that contains $p$ is maximal.  A point is called \emph{singular} if it is not  regular.

In \cite{Alex3} I introduced the concept of s.r.f.s. as the natural candidate to generalize the theory of isoparametric submanifolds in Riemannian manifolds.

\begin{definicao}[s.r.f.s.]
 Let $\F$ be a singular Riemannian foliation on a complete Riemannian manifold $M.$
$\F$ is said to be a \emph{singular Riemannian foliation with section} (s.r.f.s. for short) if for every regular point $p,$ the set $\Sigma :=\exp_{p}(\nu_p L_{p})$ is a complete immersed submanifold that meets each leaf orthogonally. $\Sigma$ is called a \emph{section}. We also suppose that the set of regular points is open and dense in $\Sigma.$ 
\end{definicao}
In Section 4 we will recall some results of \cite{Alex2} about this new theory. Among them, there is the Slice Theorem (see Theorem \ref{sliceteorema}) that claims that the restriction of s.r.f.s. to a slice is diffeomorphic to an isoparametric foliation on an Euclidean space. 
Another important result is Theorem \ref{frss-eh-equifocal} asserting that the leaves of a s.r.f.s. are locally equifocal. These two results establish the relations between the old and new theory and  allow us to associate reflections pseudogroups to a s.r.f.s. 


\section{Isoparametric and Transnormal Maps}

In the last section we present s.r.f.s. as the natural candidate to generalize isoparametric foliations on Euclidean spaces. Nevertheless, there exists another possible way to try to generalize them. We can consider the preimages of  special maps, the so called \emph{isoparametric} maps. In fact,  Cartan \cite{Cartan1,Cartan2,Cartan3,Cartan4}, Harle \cite{Harle}, Terng \cite{Terng} and Q.M. Wang \cite{Wang1} have chosen this approach.

In this section we review the definition of isoparametric and transnormal map. We also review some results about this map. Then I present  Theorem \ref{TransormalMapTeo}
that relates analytic transnormal map to s.r.f.s.  This result  can be found in \cite{Alex1}.   

 \begin{definicao}
\emph{
Let  $M^{n+q}$  be a complete Riemannian manifold. A smooth map   $H=(h_{1}\cdots h_{q}):M^{n+q}\rightarrow \mathbf{R}^{q}$ is called  \emph{transnormal} if
\begin{enumerate}
\item[0)]  $H$ has a regular value,
\item[1)]  for each regular value $c$ there exist a neighborhood $V$ of $H^{-1}(c)$ in $M$ and smooth functions $b_{i\,j}$ on $H(V)$ such that  
$<\grad h_{i}(x),\grad h_{j}(x)>= b_{i\,j}\circ H(x),$ for every $x\in V,$
\item[2)]  there is a sufficiently small neighborhood of each regular level set such that  $[\grad h_{i},\grad h_{j}]$ is a linear combination of $\grad h_{1},\cdots ,\grad h_{q},$ with coefficients being functions of $H,$ for all $i$ and $j$. 
\end{enumerate}
}
\end{definicao}
This definition is equivalent to saying that $H$ has a regular value and for each regular value $c$ there exists a neighborhood  $V$ of $ H^{-1}(c)$ in $M$ such that $H \mid_{V}\rightarrow H(V)$ is an integrable Riemannian submersion, where the Riemannian metric $g_{i\,j}$ of  $H(V)$ is the inverse matrix of $[b_{i\,j}].$

In particular, a transnormal map $H$ is said to be an \emph{isoparametric map} if  $V$ can be chosen to be $M$  and   $\bigtriangleup h_{i}= a_{i}\circ H,$ where $a_{i}$ is a smooth function.

Analytic transnormal maps can be used to describe 
 isoparametric submanifolds in a space form and equifocal submanifolds in a compact symmetric space, as we see below.

It is known (see  the book of Palais and  Terng \cite{PTlivro}) that 
given an isoparametric map $H$ on a space form and a regular value $c$ then the regular level set $H^{-1}(c)$ is an isoparametric submanifold. In addition, given an isoparametric submanifold $L$ of  $\mathbf{R}^{n+k},$ one can construct a polynomial isoparametric map from $\mathbf{R}^{n+k}$ into $\mathbf{R}^{k}$ that has $L$ as a level submanifold.

It follows from a Chevalley-type restriction theorem of E. Heintze, X. Liu and C. Olmos \cite{HOL} that 
 given an equifocal submanifold with flat sections in a simply connected symmetric space of compact type, one can find an analytic transnormal map  that has this equifocal submanifold as a regular level set.

Now we are able to state our result precisely.

\begin{teo}
\label{TransormalMapTeo}
Let $H: M\rightarrow \mathbf{R}^{q}$ be an analytic transnormal map on a real analytic complete Riemannian manifold $M.$ Let $c$ be a regular value and $L\subset H^{-1}(c)$ be a connected component of $H^{-1}(c).$ Denote  by $\Xi$  the set of all parallel normal fields along $L.$ Then 
\begin{enumerate}
\item[1)] $\F_{c,L}:=\{\eta_{\xi}(L)\}_{\xi\in\Xi}$ is a singular Riemannian foliation with sections (s.r.f.s.)  whose leaves are always embedded.
\item[2)] For each regular value $\hat{c}$ the connected components of $H^{-1}(\hat{c})$ are equifocal manifolds and leaves of $\F_{c,L}.$
\item[3)] $\F_{c,L}$ is independent of the choice of $c$ and $L$ i.e., for another regular value $\tilde{c}$ and connected component $\tilde{L}\subset H^{-1}(\tilde{c})$ we have $\F_{c,L}=\F_{\tilde{c},\tilde{L}}.$
\end{enumerate}
\end{teo}   

 Our result is  a generalization of Q.M. Wang's  result \cite{Wang1} for the case of analytic maps. Indeed, Q.M. Wang studied the so-called transnormal functions, which are real valued functions on a Riemannian manifold $M$  satisfying the equation $\|\grad f\|^{2} =b\circ f$ where $b$ is $C^{2}(f(M)).$ He  proved that the level sets of $f$ are leaves of a singular Riemannian foliation. It also follows from his proof that the regular level sets are  equifocal.  

\section{Facts about S.R.F.S.}

In this section we review  some results about s.r.f.s. that  can be found in \cite{Alex2}.

Let us start by fixing some notations.
Let $\Sigma$ be a section and $p\in\Sigma.$ Let $B$ a normal neighborhood of $p$ in $M.$ A connected component of $B\cap \Sigma$ that contains $p$ is called a \emph{local section}.

Let $p$ be a point of a leaf $L$ of a singular Riemannian foliation and let $P$ be a relatively compact, open connected neighborhood of $p$ in $L$. Then there is $\eps>0$ such that $\exp$ restricted to $\nu^\eps P=\{X\in \nu P\mid \|X\|<\eps\}$ is a diffeomorphism onto the tube $\tub (P)$ of radius $\eps$. We write $S_q=\exp(\nu^\eps_qP)$ for $q\in P$. 
We call $S_{q}$ a \emph{slice} at $q.$

A relation between a slice $S_{q}$ and  local sections  is given in the next result. 

\begin{prop}
Let $\F$ be a s.r.f.s. on a complete Riemannian manifold $M$ and let  $q\in M.$  Then 
\begin{enumerate}
\item[a)] Let $\epsilon$ be the radio of the slice $S_{q}$ and $\Lambda(q)$  the set of  local sections $\sigma$ that contain $q,$ such that  $\dist(p,q)<\epsilon$ for each $p\in\sigma.$  Then  $S_{q}= \cup_{\sigma\in\Lambda (q)}\, \sigma.$
\item[b)] $S_{x}\subset S_{q}$ for all $x\in S_{q}.$
\end{enumerate}
\end{prop}

Note that, if $q$ is a singular point, the restriction $\F|S_q$ of $\F$ to $S_{q}$ is also a singular foliation. In fact, since $\singularF$ acts transitively on the leaves,  the plaques (the connected components of the leaves intersected with $\tub(P)$) are transversal to $S_q.$

In \cite{Molino} (page 209) Molino showed that given a singular Riemannian foliation  $\F$ it is possible to change the metric  so that the restriction of $\F$ to a slice is a singular Riemannian foliation with respect to this new metric. This new metric  respects the distance between the leaves. As we see below the metric need not be changed if  the singular Riemannian foliation has sections.

\begin{cor}
Let  $\F$ be a s.r.f.s. on a complete Riemannian manifold  $M$ and  $S$  a slice. Then  $\F\cap S$  is s.r.f.s. on $ S$ with the induced metric of $M.$
\end{cor}

\begin{prop}
\label{pontos-singulares-na-geodesica}
Let $\F$ be a s.r.f.s. on a complete Riemannian manifold $M$ and
$\gamma$ be a geodesic orthogonal to a regular leaf. Then the set of singular points is isolated on $\gamma.$ 
\end{prop}

The above results were used in the proof of the following result, which is one of the main results of \cite{Alex2}.

\begin{teo}
\label{frss-eh-equifocal}
Let $\F$  be a s.r.f.s. on a complete Riemannian manifold  $M.$ Then the regular leaves are locally equifocal. In particular, the union of  the regular leaves that are equifocal is an open and dense set  in $M,$ provided that all the leaves are compact.  
\end{teo}

A consequence of this theorem is that given a regular leaf $L$ with trivial holonomy, it is possible to   reconstruct $\F$ by taking all parallel submanifolds of $L.$ More precisely we have

\begin{cor}
\label{cor-map-paralelo}
Let $L$ be a regular leaf of  $\F.$ 
\begin{enumerate}
\item[a)] Let  $\beta$ be a smooth curve of $L$ and  $\xi$ a parallel normal field along  $\beta$. Then the curve $\eta_{\xi}\circ\beta$ belongs to a leaf of \F. 
\item[b)] Let  $L$ be a regular leaf with trivial holonomy and $\Xi$ denote the set of all parallel normal fields along $L.$ Then $\F=\{\eta_{\xi}(L)\}_{\xi\in \, \Xi}.$ 
In particular, if $\xi$ is a parallel normal field along $L$ then the end point map $\eta_{\xi}:L\rightarrow L_{q}$ is surjective, where $q=\eta_{\xi}(x)$ for $x\in L.$
\end{enumerate}
\end{cor}

Corollary \ref{cor-map-paralelo} allows us to define singular holonomy maps (see definition below) which will be very useful to study $\F.$ 

Theorem \ref{frss-eh-equifocal} is also used to prove the following result, which is the other main result in \cite{Alex2}.

\begin{teo}[ Slice Theorem ]
\label{sliceteorema}
Let $\F$ be a s.r.f.s. on a complete Riemannian manifold $M$  and  $S_{q}$ the slice of a point  $q\in M.$ Then  $\F$ restrict to  $S_{q}$ is diffeomorphic to an isoparametric foliation on an open set of $\mathbf{R}^{n},$  where $n$ is the dimension of   $S_{q}.$ 
\end{teo}

\begin{cor}
\label{estratificacao-singular}
Let $\F$ be a s.r.f.s. on a complete Riemannian manifold $M$ and $\sigma$ be a local section. Then the set of singular points of $\F$ contained in $\sigma$ is a finite union of totally geodesic hypersurfaces that are diffeomorphic to  focal hyperplanes contained in a section of an  isoparametric foliation on an open set of an euclidian space.   
\end{cor}
This set of singular points of $\F$ contained in $\sigma$ is going to be called 
 \emph{singular stratification of the local section.}

Denote by $M_{r}$ the set of regular points in $M.$ A \emph{Weyl Chamber} of a local section $\sigma$ is the union of a connected component $A$  of $M_{r}\cap \sigma$ with the walls of the singular stratum of $\sigma$ that are contained in $\partial A$.

\begin{cor}
\label{convexidadedeWeyl}
A Weyl Chamber of a local section is a convex set.
\end{cor}

Another consequence of the Slice Theorem is the following result, which is a converse of Theorem \ref{TransormalMapTeo}.

\begin{cor}[Trivialization of $\F$]
\label{frss-transnormal}
The plaques of a s.r.f.s. are always level sets of a transnormal map.
\end{cor}

The slice theorem gives us a description of  the plaques of a s.r.f.s. However, it does not guarantee that two different plaques belong to  different leaves.
To obtain this  kind of information, we must extend the concept of holonomy to describe not only what happens near to a regular leaf, but also what happens in a neighborhood of a singular leaf.

\begin{prop}[Singular Holonomy]
\label{prop-holonomia-singular}
Let $L_{p}$ be a regular leaf,  $\beta$ a smooth curve in $L_{p}$ and let $[\beta]$ denote the homotopy class of $\beta.$ Let $U$ be a local section that contains $p=\beta(0).$ Then there exists a local section $V$ that contains $\beta(1)$ and 
an isometry $\varphi_{[\beta]}:U\rightarrow V$  that has the following properties:
\begin{enumerate}
\item[1)]$\varphi_{[\beta]}(x)\in L_{x}$ for each $x\in U,$
\item[2)]$d\varphi_{[\beta]}\xi(0)=\xi(1),$ where $\xi(s)$ is a  parallel normal field along $\beta(t).$
\end{enumerate}
\end{prop}

An isometry as above is called \emph{singular holonomy map along $\beta$}. 
 We remark that, in the definition of the singular holonomy map, the domain $U$ can contain singular points. If the domain $U$ and the range $V$ are sufficiently small, then the singular holonomy map  coincides with the holonomy map along $\beta.$

After the Slice Theorem and its corollaries, one can ask if there are reflections in the walls of the Weyl Chambers or something that is similar to the  Weyl groups of  the classical Lie Group theory. The next definition and results answer this question.

\begin{definicao}[Weyl Pseudogroup $W$]
\label{definitionWeylPseudogroup}
 The pseudosubgroup  generated by all singular holonomy maps $\varphi_{[\beta]}$ such that $\beta(0)$ and $\beta(1)$ belong to the same local section $\sigma$ is called \emph{generalized Weyl pseudogroup of} $\sigma.$ Let $W_{\sigma}$ denote this pseudogroup. Similarly we define $W_{\Sigma}$ for a section $\Sigma.$ 
\end{definicao}

\begin{rem}
To recall the definition of pseudogroups and  orbifolds  see E. Salem \cite[Appendix D]{Molino}. 
\end{rem}

\begin{prop}
Let $\F$ be a s.r.f.s. on a complete Riemannian manifold $M.$
Let $\sigma$ be a local section. Then the reflections in the hypersurfaces of the singular stratum of the local section $\sigma$ leave $\F|\sigma$  invariant. Moreover these reflections are elements of $W_{\sigma}.$
\end{prop}

One can construct examples of s.r.f.s by suspension such that  $W_{\sigma}$ is larger than the pseudogroup generated by the reflections  in the hypersurfaces of the singular stratum of $\sigma.$ The next result gives a sufficient condition to guarantee that both pseudogroups coincide.

\begin{prop}
\label{Weylgrouparereflection}
Let $\F$ be a s.r.f.s. on a complete Riemannian manifold $M.$
Suppose that each leaf is compact has trivial holonomy. Let $\sigma$ be a local section.  Then  $W_{\sigma}$ is generated by the reflections  in the hypersurfaces of the singular stratum of the local section. 
\end{prop}

\begin{cor}
 Let $\F$ be a s.r.f.s. on a complete Riemannian manifold $M.$ Suppose that the leaves are compact and  the holonomies of regular leaves are trivial. 
Consider a tubular neighborhood $\tub(L_{q})$ of a leaf $L_{q}$ and let $L_{p}$ be a regular leaf contained in $\tub(L_{q}).$ Finally denote by 
 $\Pi:\tub(L_{q})\rightarrow L_{q}$ the orthogonal projection. Then $L_{p}$ is the total space of a fiber bundle with a projection  $\Pi,$ a basis   $L_{q}$ and a fiber diffeomorphic to an isoparametric submanifold of an Euclidean space.
\end{cor}

\begin{prop}
\label{estrutura-transversa-do-fecho} 

Let $\F$ be a s.r.f.s. on a complete Riemannian manifold $M$ and let $\sigma$ be a local section. Consider a point  $p\in\sigma,$ then
\[\overline{W_{\sigma}}\cdot p= \overline{L_{p}\cap\sigma}.\]
In other words, the closure of $L_{p}\cap\sigma$ is an orbit of 
complete closed pseudogroup of local isometries. In particular $\overline{L_{p}\cap\sigma}$ is a closed submanifold.
\end{prop}

\section{Exceptional Leaves and the Space of Orbits of a S.R.F.S}


Proposition \ref{Weylgrouparereflection} motivates us to look for sufficient conditions for  the holonomy of each regular leaf of s.r.f.s. to be trivial. An answer to this problem is suggested by the Lie Group theory. In fact, it is known that a simply connected compact Lie group $G$  has no exceptional conjugacy class (see Duistermaat and Kolk \cite[pag 164, 3.9.5 Corollary]{Duistermaat}). 
Therefore it is natural to ask if the holonomy of regular leaves of a s.r.f.s. on $M$ is trivial  if $M$ is simply connected.


In \cite{AlexToeben} T\"{o}ben and I proved the following result:

\begin{teo}
\label{teo-excepcional-leaf}
Let $\F$ be a s.r.f.s. on a simply connected complete Riemannian manifold $M.$ Suppose that the leaves of $\F$ are compact. Then each regular leaf has trivial normal holonomy.
\end{teo}

This result generalizes a theorem of  T\"{o}ben \cite{Toeben2}, who proved the result under the additional assumption that the sections are symmetric spaces or do not have any conjugate points.

\begin{rem}
Note that if $M$ is not simply connected space, then the leaves can also have trivial holonomy. For example consider $M=\mathbf{S}^{1}\times\mathbf{S}^{1}$ and define $\F$ as  $\{\mathbf{S}^{1}\times p\}_{p\in\mathbf{S}^{1}}.$  
\end{rem}

Now consider the  adjoint action of $\mathbf{SU(2)}.$ By direct calculations (see \cite[pag 11]{Duistermaat}), it is easy to see that, for each section $\Sigma$, the closure of a connected component $\Omega$ of the set of regular points in $\Sigma$ is a fundamental domain, i.e., each orbit of the adjoint action intersect $\overline{\Omega}$ exactly  one time. 
This is a special case of the following joint result  with T\"{o}ben \cite{AlexToeben}.





\begin{teo}
\label{teo-fundamental-domain}
Let $\F$ be a s.r.f.s. on a simply connected Riemannian manifold $M.$ Suppose also that the leaves of $\F$ are compact. Then
\begin{enumerate}
\item $M/F$  is a simply connected Coxeter orbifold.
\item Let $\Sigma$ be a section of $\F$ and $\Pi:M\rightarrow M/\F$ the projection. Let $\Omega$ be a connected component of the set of regular points in $\Sigma.$ Then  $\Pi:\Omega\rightarrow M_{r}/\F$ and  $\Pi:{\overline{\Omega}}\rightarrow M/\F$ are homeomorphisms. In addition $\Omega$ is convex, i.e., for each $p,q$  in $\Omega$  each minimal segment of geodesic that joins $p$ to $q$ is contained in $\Omega.$
\end{enumerate}
\end{teo}

\section{Equifocal Submanifolds and S.R.F.S}

 Theorem \ref{frss-eh-equifocal} asserts that the leaves of a s.r.f.s. are local equifocal.
It is natural to ask if the converse is true, i.e., if the parallel submanifolds of an equifocal submanifold are leaves of a s.r.f.s. 

As we have said before, Terng and Thorbergsson \cite{TTh1}
answered this question when the equifocal submanifold has flat section and is embedded in a symmetric space of compact type.
I also answered this question in \cite{Alex1} when the equifocal submanifold is image of a regular value of an analytic transnormal map (see Theorem \ref{TransormalMapTeo}). 
Recently T\"{o}ben  \cite{Toeben1} used the blow up technique to study local equifocal submanifolds (which he called submanifold with parallel focal structure). 
He gave a necessary and sufficient condition for a closed embedded local equifocal submanifold to induce a s.r.f.s. 


In this section I will give an alternative proof of   T\"{o}ben's Theorem. This proof relies on  \cite{Alex1}.  

To begin with, let us recall some definitions.


As already said, a local section $\sigma$ contained in a section $\Sigma$ is a connected component of $B\cap \Sigma,$ where $B$ is a normal neighborhood in $M.$ We say that two local sections $\sigma$ and $\tilde{\sigma}$ are the same if there exists a $p\in \sigma\cap\tilde{\sigma}$ and $T_{p}\sigma=T_{p}\tilde{\sigma}.$

Let $\Sigma$ be a section of $L$ and $p$ a point of $\Sigma.$ The point $p$ is called \emph{regular} if there exists only one section in $\Sigma$ that contains $p,$ i.e., if given local sections $\sigma$ and $\tilde{\sigma}$  that contain $p,$ then $\sigma=\tilde{\sigma}.$

Now we can state the T\"{o}ben's result (see \cite[pag 4,Theorem A]{Toeben1}) for the case of equifocal submanifolds.
\begin{teo}[T\"{o}ben]
A closed and embedded equifocal submanifold $L$ is a leaf of a s.r.f.s. $\F$ if and only there is only one section through each regular point of the normal exponential map of $L.$ 
\end{teo}

As T\"{o}ben remarked, one direction of the theorem (the necessary condition) was already proved in \cite{Alex2}. In fact, if $\F$ is a s.r.f.s. then it follows from Theorem \ref{frss-eh-equifocal}
that $L$ is an local equifocal submanifold. The fact that there is only one section through each regular point  of the normal exponential map of $L$ also follows from Theorem \ref{frss-eh-equifocal} (see details in \cite{Alex2}).

 On the other hand, as T\"{o}ben also remarked (\cite[pag 19]{Toeben1}), the set of regular points of the normal exponential map is an open and dense set in each section. Hence the condition that there is only one section through each regular point of normal exponential map  implies that the set of regular points of each section $\Sigma$  (in the sense  defined above) is an open and dense set in $\Sigma.$ Therefore all we have to prove is the following result. 

\begin{teo}
\label{teo-equical-srfs}
Let $L$ be a closed embedded equifocal submanifold.  Suppose that set of regular points of each section $\Sigma$ is an open and dense set in $\Sigma.$ Define  $\Xi$ as the set of all parallel normal fields along $L.$ Then $\F:=\{\eta_{\xi}(L)\}_{\xi\in \, \Xi}$ is a s.r.f.s.
\end{teo}

In order to prove this result, we need the following results.

\begin{prop}[Heintze, Olmos and Liu \cite{HOL}]
\label{lemma-HOL}
Let $M$ be a complete Riemannian manifold, $L$ be an immersed submanifold of  $M$ with  globally flat normal bundle and  $\xi$ be a normal parallel field along $L.$ Suppose that  $\Sigma_{x}:=\exp_{x}(\nu_{x}L)$ is a totally geodesic complete submanifold for all $x\in L,$ i.e., $L$ has sections. Then 
\begin{enumerate}
\item $ d \eta_{\xi}(v)$ is orthogonal to $\Sigma_{x}$ at $\eta_{\xi}(x)$ for all $v\in T_{x}L.$
\item Suppose that $\exp_{p}(\xi)$ is not a focal point of $p\in L.$ Then there exists  a neighborhood $U$ of $p$ in $L$ such that  $\eta_{\xi}(U)$ is an embedded submanifold, which meets $\Sigma_{x}$ orthogonally and has globally flat normal bundle. In addition,  a parallel normal field along $U$ transported to $\eta_{\xi}(U)$ by parallel translation along the geodesics $\exp(t \xi)$ is a parallel normal field along $\eta_{\xi}(U).$
\end{enumerate}
\end{prop} 

\begin{prop}
\label{prop-slice}
Let $L$ be a closed embedded equifocal  submanifold in a complete Riemannian manifold $M.$ Suppose that set of regular points of each section $\Sigma$ is an open and dense set in $\Sigma.$ Let $p$ be a point in $L$ and let $U$ be a neighborhood of $p$ in $L$ such that $\eta_{\xi}(U)$ is an embedded submanifold. Set $q:=\eta_{\xi}(p).$ Consider $\tub(\eta_{\xi}(U))$ a tubular neighborhood of $\eta_{\xi}(U).$   
 Let  $S_{q}$ denote the connected component of  $\exp_{q}(\nu_{q}\eta_{\xi}(U))\cap \tub (\eta_{\xi}(U))$ that contains $q.$  Then  
\[S_{q}= \cup_{\sigma\in\Lambda(q)} \, \sigma\]
where  $\Lambda(q)$ is the set of local sections that contain  $q.$ 
\end{prop}

\begin{proof} 
Let $F$ be a connected component of $\eta_{\xi}^{-1}(q)$ that contains p. Choose a constant $0<k<1$ and define $\xi_{1}:=k\,\xi$ such that
\begin{enumerate}
\item $\eta_{\xi_{1}}(p)$ belongs to $\tub(\eta_{\xi}(U)).$ 
\item $\dim\eta_{\xi_{1}}(U)=\dim U.$ 
\end{enumerate}

Let  $F_{1}:=\eta_{\xi_{1}}(F).$ Note that $F_{1}$ is an embedded compact submanifold contained in $\tub(\eta_{\xi}(U)).$

The following claim can easily be checked.

\begin{claim}
\label{claim1}
For each $z\in F$ there exists a neighborhood $\widetilde{U}$ of $z$ in $L$ such that $\eta_{\xi}(\widetilde{U})$ is a neighborhood of $\eta_{\xi}(U).$ 
\end{claim}

\begin{claim}
\label{claim2}
Let $z_{1}$ be a point of $F_{1}$ and $\sigma$ be a local section that contains $z_{1}.$ Then $\sigma\subset S_{q}$ and $q\in\sigma.$ In particular $F_{1}\subset S_{q}.$
\end{claim}
\begin{proof}
Let $z\in F$ such that $\eta_{\xi_{1}}(z)=z_{1}.$ 
It follows from Proposition \ref{lemma-HOL} that, for a small neighborhood $\widetilde{U}$ of $z$ in $L,$  $\eta_{\xi}(\widetilde{U})$ is orthogonal to the section $\Sigma$ that contains $z.$ (Note that $\Sigma$ also contains $z_{1}$).

On the other hand, it follows from Claim \ref{claim1} that $\eta_{\xi}(\widetilde{U})$ is a neighborhood of $\eta_{\xi}(U)$ for a small neighborhood $\widetilde{U}$ of $z.$ Therefore, $\Sigma$ is orthogonal to $\eta_{\xi}(U)$ and this implies the claim.
\end{proof}

\begin{claim}
\label{claim3}
Let $z_{1}\in F_{1}$ and $\sigma$ be a local section that contains $z.$ Then 
\[T_{z_{1}}S_{q}=T_{z_{1}}\sigma  \oplus   T_{z_{1}}F_{1}\]
\end{claim}
\begin{proof}
Let $z\in F$ such that $\eta_{\xi_{1}}(z)=z_{1}.$ Let $\Sigma$ be the section that contains $z$ and hence that contains $\sigma.$ Let $\tilde{U}$ be a neighborhood of $z$ in $L$ as defined in Claim \ref{claim1}. 
It is easy to see that
\begin{eqnarray*}
\dim \Sigma +\dim F +\dim \eta_{\xi}(\widetilde{U})&=& \dim M\\
                                                   &=&\dim S_{q}+ \dim \eta_{\xi}(\widetilde{U}).
\end{eqnarray*}
Since $\dim\sigma=\dim\Sigma$ and $\dim F=\dim F_{1}$ we have
\[ \dim \sigma+\dim F_{1}=\dim S_{q}.\]
Now our claim follows from Claim \ref{claim2} and from the above equation.
\end{proof}

\begin{claim}
\label{claim4}
$S_{q}\subset \cup_{\sigma\in\Lambda(q)} \, \sigma$
\end{claim}
\begin{proof}
Let $p\in S_{q}$ and let $\gamma$ be the shortest segment of geodesic in $S_{q}$ that joins $p$ to $F_{1}$ with $\gamma(0)=p$ and $\gamma(1)\in F_{1}.$ (Note that $\gamma$ is a geodesic in $S_{q}$ but  it might not be a geodesic in $M$). Since $\gamma$ is orthogonal to $F_{1},$  
 Claim \ref{claim2} and Claim \ref{claim3} imply that $\gamma\subset\sigma,$ where $\sigma$ is the local section that contains $\gamma(1).$ The fact that $\gamma\subset \sigma$ implies that  $p\in\sigma.$ 

On the other hand, it follows from Claim \ref{claim2} that $q\in\sigma.$ Therefore $p\in\sigma\in \Lambda(q).$  
\end{proof}

\begin{claim}
\label{claim5}
$S_{q}\supset \cup_{\sigma\in\Lambda(q)} \, \sigma$
\end{claim}
\begin{proof}
Let $\sigma\in\Lambda(q).$ Then there exists a regular point $p\in\sigma.$ Let $\tilde{q}\in\eta_{\xi}(U)$ such that $p\in S_{\tilde{q}}.$ It follows from the proof of Claim \ref{claim4} that there exists a local section $\tilde{\sigma}\subset S_{\tilde{q}}$ such that $p\in\tilde{\sigma}$ and $\tilde{q}\in\tilde{\sigma}.$ Since $p$ is regular, $\tilde{\sigma}=\sigma$ and hence $q\in S_{\tilde{q}}.$ Since $\tub(\eta_{\xi}(U))$ is a tubular neighborhood, this is possible only if $S_{q}=S_{\tilde{q}}.$ Therefore $\sigma=\tilde{\sigma}\subset S_{\tilde{q}}=S_{q}.$
\end{proof}

\end{proof}

\begin{cor}
\label{cor-propslice}
 Let $\sigma$ be a local section such that $\sigma\cap S_{q}\neq \emptyset.$ Then $\sigma\subset S_{q}$ and $q\in \sigma.$
\end{cor} 
\begin{proof}
The proof goes exactly in the proof of Claim \ref{claim5}.
\end{proof}


\textbf{Proof of Theorem \ref{teo-equical-srfs}} 

Now the proof of Theorem \ref{teo-equical-srfs} is basically the same proof of Item 1 of Theorem \ref{TransormalMapTeo} apart from the modified Proposition \ref{prop-slice}. Therefore we only sketch its main steps. For details (e.g. the proofs of Corollary \ref{lema-transnormalmap-frss}, Claim \ref{claimproof1}, Claim \ref{claimproof2} and Equation \ref{SxsubsetSy}) see \cite{Alex1}.

Proposition \ref{prop-slice} implies the following corollary.

\begin{cor}
\label{lema-transnormalmap-frss}
Let   $U\subset L$ and   $\tilde{U}\subset L$ denote two open sets small enough such that $\eta_{\xi}(U)$ and $\eta_{\tilde{\xi}}(\tilde{U})$ are embedded.
Let $\xi$ and $\tilde{\xi}$ be normal parallel fields along  $U$ and $\tilde{U}$ respectively.  Suppose that 
\begin{enumerate}
\item[i)]$\eta_{\tilde{\xi}}(\tilde{U})\subset \tub (\eta_{\xi}(U)),$
\item[ii)]  $\eta_{\tilde{\xi}}(\tilde{U})$ is not contained in  $\eta_{\xi}(U).$
\end{enumerate}
Then $\eta_{\tilde{\xi}}(\tilde{U})$ is equidistant from  $\eta_{\xi}(U),$ meaning that it is contained in a boundary of a tube  with axis $\eta_{\xi}(U).$
\end{cor}
\begin{rem}
\emph{ Note that we do not assume that $\eta_{\tilde{\xi}}(\tilde{U})$ and $\eta_{\xi}(U)$ are disjoint, this will be a consequence of the corollary.}
\end{rem}

It follows from  Corollary \ref{lema-transnormalmap-frss}  that  $\eta_{\xi}:L\rightarrow M$ is an immersion without self intersection. Since it is also proper, we conclude that 
$\eta_{\xi}(L)$ is embedded. 

Corollary  \ref{lema-transnormalmap-frss} implies that $\F:=\{\eta_{\xi}(L)\}_{\xi\in \, \Xi}$ is a partition of  $M.$ 

Proposition \ref{lemma-HOL} implies that the sections of $L$ are orthogonal to the leaves of $\F.$  

Now we have to check that the partition $\F$ is a transnormal system,  i.e., every geodesic that is perpendicular at one point to a leaf remains perpendicular to every leaf it meets. Let $\eta_{\xi}(U)$ be a plaque and $\tub(\eta_{\xi}(U))$ be a tube that has  this plaque as the axis.  Then Proposition \ref{prop-slice} implies the following claims.

\begin{claim}
\label{claimproof1} 
If a segment of geodesic is perpendicular to the axis of the tube, then the geodesic is perpendicular to each plaque in the tube that it meets. 
\end{claim}

\begin{claim}
\label{claimproof2} 
If a segment of geodesic  joins a plaque to the axis and is perpendicular to this plaque, then it is perpendicular to the axis.   
\end{claim}

Now given a geodesic that is perpendicular to a leaf, we can cover it with tubular neighborhoods and we can use Claim \ref{claimproof1} and Claim \ref{claimproof2} to see that this geodesic is always perpendicular to each leaf that it meets.  

To conclude, we have only to check that $\F$ is singular, i.e., for each $p$ we have to find vector fields $\{v_{i}\}$ that are always tangent to the leaves of $\F$ and such that 
$\{v_{i}(p)\}$ is a basis of $T_{p}L.$

Let $P\subset L_{p}$  be a plaque. We can find an orthogonal frame $\{v_{i}\}$ of $\nu S_{y}$ for each $y\in P$ such that their norms are 1 near $p$ and $0$ outside of 
$\tub(P).$  

It is not difficult to conclude from Proposition \ref{prop-slice}  that 
\begin{eqnarray}
\label{SxsubsetSy}
S_{x}\subset S_{y}
\end{eqnarray}
for $x\in S_{y}.$
Since
\[S_{x}\subset S_{y}\Rightarrow \nu S_{y}\subset\nu S_{x},\]
 the fields $\{v_{i}\}$ have the desired property.

\bibliographystyle{amsplain}

\end{document}